\documentclass[12pt]{amsart}
\usepackage{amssymb}
\usepackage{verbatim}
\usepackage[usenames]{color}
\usepackage{hyperref}

\newtheorem{thm}{Theorem}
 
\newtheorem{la}[thm]{Lemma}
\newtheorem{cor}[thm]{Corollary}

\theoremstyle{definition}
\newtheorem{df}[thm]{Definition}
\newtheorem{conv}[thm]{Convention}

\theoremstyle{remark}

\DeclareMathOperator{\nfrk}{nfrk}

\newenvironment{eqnum}{\begin{equation}}{\end{equation}}

\newenvironment{ls}{\begin{itemize}}{\end{itemize}}
\newenvironment{lsnum}{\begin{enumerate}}{\end{enumerate}}
\newenvironment{pf}{\begin{proof}}{\end{proof}}

\newcommand{\scr}[1]{\ensuremath{\mathcal {#1}}}

\newcommand{\bbb}[1]{\ensuremath{\mathbb {#1}}}

\renewcommand{\phi}{\varphi}

\newcommand{\sq}[1]{\ensuremath{\langle#1\rangle}}
\newcommand{\fil}[1]{\ensuremath{(#1)}}

\newcommand{\restr}{\mathop{\upharpoonright}}

\newcommand{\notarrow}{\kern .42em\not\kern -.42em\longrightarrow}
\renewcommand{\th}{\ensuremath{{}^{\text{th}}}}
\newcommand{\nf}[1]{\ensuremath{\nfrk(#1)}}
\newcommand{\zp}{\ensuremath{\bbb Z_p}}
\newcommand{\qp}{\ensuremath{\bbb Q_p}}

\newcommand{\noprint}[1]{\relax}

\title{Disjoint Non-Free Subgoups of Abelian Groups}
\author{Andreas Blass}
\thanks{Blass was partially supported by the United States National
  Science Foundation, grant DMS--0070723}
\address{Mathematics Department\\
University of Michigan\\
Ann Arbor, MI 48109--1043, U.S.A.}
\email{ablass@umich.edu}
\author{Saharon Shelah}
\address{Mathematics Department\\
Hebrew University\\
Jerusalem 91904, Israel\\
and
Mathematics Department\\
Rutgers University\\
New Brunswick, NJ 08903, U.S.A.}
\thanks{Shelah was partially supported by the US-Israel Binational
  Science Foundation, grant 2002323.  This is paper 870 in his list of
publications}

\begin{document}

\begin{abstract}
Let $G$ be an abelian group and let $\lambda$ be the smallest
rank of any group whose direct sum with a free group is isomorphic to
$G$.  If $\lambda$ is uncountable, then $G$ has $\lambda$ pairwise
disjoint, non-free subgroups.  There is an example where $\lambda$ is
countably infinite and $G$ does not have even two disjoint, non-free
subgroups.  
\end{abstract}

\maketitle

\section{Introduction}

In a discussion between the first author and John Irwin, the question
arose whether every non-free, separable, torsion-free abelian group
has two disjoint non-free subgroups.  (Of course, in the context of
subgroups, ``disjoint'' means that the intersection is $\{0\}$.)  The
main result of this paper is a strong affirmative answer.  To state
the result in appropriate generality, we need some terminology.

\begin{conv}
  All groups in this paper are understood to be abelian and
  torsion-free.  In particular, ``free group'' means ``free abelian
  group''.
\end{conv}

\begin{df}
  The \emph{non-free rank} of a group $G$, written \nf G, is the
  smallest cardinal $\kappa$ such that $G$ can be split as the direct
  sum of a group of rank $\leq\kappa$ and a free group.
\end{df}

\begin{thm}   \label{main}
  If \nf G is uncountable, then $G$ has \nf G pairwise disjoint,
  non-free subgroups.
\end{thm}

Recall that any countable, separable group is free.  It follows that,
if $G$ is separable and not free, then \nf G is necessarily
uncountable, so the theorem applies to $G$.  It gives not only two
disjoint non-free subgroups as in the original question but $\nf
G\geq\aleph_1$ of them.

The number of disjoint, non-free subgroups obtained in the theorem is
the most one could hope for.  Indeed, if $G\cong H\oplus F$ where $F$
is free and $H$ has rank and therefore cardinality equal to the
infinite cardinal \nf G, then any non-free subgroup $S$ of $G$ must
have a non-zero intersection with $H$.  Otherwise the projection to
$F$ would map $S$ one-to-one into the free group $F$ and it would
follow that $S$ is free.  Therefore, disjoint non-free subgroups of
$G$ must intersect $H-\{0\}$ in disjoint non-empty sets.  So there
cannot be more such subgroups than $|H|=\nf G$.

Although the theorem gives an optimal result for separable groups, the
fact that it does not explicitly mention separability raises another
question: Is the uncountability hypothesis really needed?  We shall
answer this question affirmatively by exhibiting a (necessarily
non-separable) non-free group $G$ such that $\nf G=\aleph_0$ and $G$
does not have two disjoint non-free subgroups, let alone $\aleph_0$ of
them.

This paper contains, in addition to this introduction and a section of
known preliminary results, three sections.  The first two are devoted
to the proof of Theorem~\ref{main}.  Section~\ref{reg} contains the
proof for the case that \nf G is (uncountable and) regular.
Section~\ref{sing} contains the additional arguments needed to extend
the result to the case of singular \nf G.  Finally,
Section~\ref{counter} presents our counterexample for the case of
countable non-free rank.

\section{Preliminaries}   \label{prelim}

In this section, we collect for reference some conventions,
definitions, and known results that will be needed in our proofs.  The
book \cite{em} of Eklof and Mekler serves as a standard reference for
this material.

\begin{conv}
Group operations will be written additively.  For $n\in\bbb Z$ and $x$
  a group element, the notation $nx$ means the sum of $n$ copies of
  $x$ if $n>0$, it means $(-n)(-x)$ if $n<0$, and it means 0 if $n=0$.
\end{conv}

\begin{df}
A subgroup $H$ of $G$ is \emph{pure} in $G$ if, whenever $x\in G$ and
$nx\in H$ for some non-zero integer $n$, then $x\in H$.  If $H$ is an
arbitrary subgroup of $G$, then we write $H_*$ for the
\emph{purification} of $H$, the smallest pure subgroup of $G$ that
includes $H$.

If $X$ is any subset of a group $G$, we write \sq X for the subgroup
of $G$ generated by $X$.  Its purification $\sq X_*$ is called the
\emph{pure subgroup generated by} $X$.

For a prime number $p$, the field of $p$-adic numbers will be denoted
by \qp, and the subring of $p$-adic integers will be denoted by \zp.
\end{df}

\begin{conv}
  When we refer to the numerators or denominators of rational numbers,
  we always mean that the rational numbers are regarded as fractions
  in reduced form.  
\end{conv}

Recall that the rational field \bbb Q is a subfield of \qp\ and that
the intersection $\bbb Q\cap\zp$ consists of those rational numbers
whose denominators are not divisible by $p$.  A useful consequence is
that a rational number is in \zp\ for all primes $p$ if and only if it
is an integer.  Recall also that $p\zp$ is the unique maximal ideal of
\zp\ and that the quotient $\zp/p\zp$ is isomorphic to the $p$-element
field $\bbb Z/p\bbb Z$, the isomorphism being induced by the inclusion
of \bbb Z in \zp.

Let $\lambda$ be an uncountable regular cardinal and $G$ a group of
cardinality $\lambda$.  A \emph{filtration} of $G$ is an increasing
sequence \fil{G_\alpha:\alpha<\lambda} of subgroups of $G$, each of
cardinality $<\lambda$, continuous at limit ordinals $\beta<\lambda$
(i.e., $\bigcup_{\alpha<\beta}G_\alpha=G_\beta$), and with
$\bigcup_{\alpha<\lambda}G_\alpha=G$.  Although there are many
filtrations of $G$, any two of them, say \fil{G_\alpha} and
\fil{G'_\alpha}, agree almost everywhere in the sense that the set of
agreement $\{\alpha<\lambda:G_\alpha=G'_\alpha\}$ includes (in fact
is) a closed unbounded set (\emph{club}) in $\lambda$.

Recall that a subset of $\lambda$ is \emph{stationary} if it
intersects every club and that the intersection of any fewer than
$\lambda$ clubs is again a club.  It follows that a stationary set
cannot be the union of fewer than $\lambda$ non-stationary sets.  We
shall also need Fodor's theorem and a variant of it.  Fodor's theorem
says that, if $S$ is a stationary subset of $\lambda$ and
$f:S\to\lambda$ is a regressive function (i.e., $f(\alpha)<\alpha$ for
all $\alpha\in S$), then $f$ is constant on a stationary subset of
$S$.  The variant that we shall need is the following.

\begin{la}   \label{multifodor}
  Let \fil{G_\alpha} be a filtration of $G$, let $S\subseteq\lambda$ be
  stationary, and let $f:S\to G$ be a function such that $f(\alpha)\in
  G_\alpha$ for all $\alpha\in S$.  Then $f$ is constant on some
  stationary subset of $S$.
\end{la}

\begin{pf}
Replacing $S$ by its intersection with the club of all limit ordinals,
we may assume that every $\alpha\in S$ is a limit ordinal.  For such
$\alpha$, continuity of the filtration tells us that $f(\alpha)\in
G_{g(\alpha)}$ for some $g(\alpha)<\alpha$.  By Fodor's theorem, the
regressive function $g$ is constant, say with value $\beta$, on some
stationary subset $T$ of $S$.  Since $|G_\beta|<\lambda$, we have a
decomposition of $T$ into fewer than $\lambda$ sets $T_x=\{\alpha\in
T: f(\alpha)=x\}$ for $x\in G_\beta$.  So one of the pieces $T_x$
must be stationary, and the lemma is established.
\end{pf}

The connection between filtrations and freeness is given by the
following definition and lemma.

\begin{df}   \label{gamma-df}
  The \emph{Gamma invariant} of $G$ is 
\[
\Gamma(G) = \{\alpha<\lambda: G/G_\alpha\text{ has a non-free subgroup
of cardinality }<\lambda\}.
\]
\end{df}

This definition seems to depend on the choice of filtration, but in
fact it doesn't, modulo restriction to a club, because any two
filtrations agree on a club.  In particular, the statement
``$\Gamma(G)$ is stationary'' has the same truth value for all choices
of the filtration.

\begin{la}   \label{gamma-la}
$\Gamma(G)$ is stationary if and only if $\nf G=\lambda$.
\end{la}

\begin{pf}
Suppose first that $\Gamma(G)$ is not stationary and is therefore
disjoint from some club $C$.  Let $c:\lambda\to C$ be the function
enumerating $C$ in increasing order.  We construct, by induction on
$\alpha$, a free basis $B_\alpha$ for $G_{c(\alpha)}/G_{c(0)}$ such
that $B_\alpha\subseteq B_\beta$ whenever $\alpha<\beta$.  This is
trivial for $\alpha=0$.  Continuity of $c$ (because $C$ is closed)
lets us simply take unions at limit stages.  At a successor step, say
from $\alpha$ to $\alpha+1$, remember that $c(\alpha)\in
C\subseteq\lambda-\Gamma(G)$, so $G_{c(\alpha+1)}/G_{c(\alpha)}$ has a
free basis.  Pick one representative in $G_{c(\alpha+1)}$ for each
member of this basis, and adjoin the chosen representatives to
$B_\alpha$ to get $B_{\alpha+1}$.  At the end of the induction,
$\bigcup_{\alpha<\lambda}B_\alpha$ is a free basis for $G/G_{c(0)}$.
Since homomorphisms onto free groups split, $G$ is isomorphic to the
direct sum of $G_{c(0)}$ and a free group.  So $\nf
G\leq|G_{c(0)}|<\lambda$.  

For the converse, suppose $G=H\oplus F$ where $|H|<\lambda$ and $F$ is
free.  Fix a free basis for $F$, necessarily of cardinality $\lambda$,
and enumerate it in a sequence of order-type $\lambda$.  Then $G$ has
a filtration whose $\alpha\th$ element $G_\alpha$ is the subgroup
generated by $H$ and the first $\alpha$ elements in the enumeration of
the basis of $F$.  Then each of the quotients $G/G_\alpha$ is free,
with a basis represented by all but the first $\alpha$ elements of the
basis of $F$.  Therefore, all subgroups of $G/G_\alpha$ are also
free.  So $\Gamma(G)$ is empty for this filtration, and therefore
non-stationary for all filtrations.
\end{pf}

It will be convenient to use filtrations normalized as follows.

\begin{la}   \label{fil-norm}
  $G$ has a filtration \fil{G_\alpha:\alpha<\lambda} such that
\begin{ls}
  \item each $G_\alpha$ is a pure subgroup of $G$, and
  \item whenever $G/G_\alpha$ has a non-free subgroup of cardinality
$<\lambda$ (i.e., whenever $\alpha\in\Gamma(G)$ as calculated with
this filtration), then $G_{\alpha+1}/G_\alpha$ is such a subgroup.
\end{ls}
\end{la}

\begin{pf}
Starting with any filtration \fil{H_\alpha} of $G$, we produce a new
filtration with the desired additional properties by defining $G_\alpha$
inductively.  Start with $G_0=(H_0)_*$, and at limit ordinals take unions
(as demanded by the definition of filtration).  At a successor step
from $\alpha$ to $\alpha+1$, first choose a subgroup $K$ of $G$, of
cardinality $<\lambda$, such that $G_\alpha\subsetneq K$ and such
that, if possible, $K/G_\alpha$ is not free.  Then let
$G_{\alpha+1}=(H_{\alpha+1}+K)_*$.  
\end{pf}

\section{Proof for Regular Non-Free Rank}  \label{reg}

In this section, we shall establish Theorem~\ref{main} in the case
that \nf G is an uncountable, regular cardinal $\lambda$.  We may
assume, without loss of generality, that $|G|=\lambda$.  Indeed, if
$G$ were larger, we could use, in place of $G$, the summand $H$ in the
decomposition $G\cong H\oplus F$ given by the definition of \nf G.
This $H$ has rank, cardinality, and non-free rank all equal to
$\lambda$, and of course if we find $\lambda$ pairwise disjoint
non-free subgroups in $H$ then these will serve in $G$ as well.  From
now on, assume $|G|=\lambda$.

For the rest of this section, we fix a filtration
\fil{G_\alpha:\alpha<\lambda} with the properties in
Lemma~\ref{fil-norm}.  We define $\Gamma(G)$ using this filtration in
Definition~\ref{gamma-df}, and we note that, by Lemma~\ref{gamma-la},
$\Gamma(G)$ is a stationary subset of $\lambda$.

For each $\alpha\in\Gamma(G)$, the properties of $(G_\alpha)$ in
Lemma~\ref{fil-norm} imply that $G_{\alpha+1}/G_\alpha $ is a
torsion-free, non-free group.  Let $Y_\alpha$ be a maximal linearly
independent subset of this group.  By linear independence, the
subgroup \sq{Y_\alpha} generated by $Y_\alpha$ is free, and by
maximality, its purification $\sq{Y_\alpha}_*$ is all of
$G_{\alpha+1}/G_\alpha $.  Choose, for each element of $Y_\alpha$, a
representative in $G_{\alpha+1}$, and let $X_\alpha$ be the set of
these chosen representatives.  Thus, the projection from
$G_{\alpha+1}$ to its quotient modulo $G_\alpha$ maps $X_\alpha$
one-to-one onto $Y_\alpha$.

Expressing the properties of $Y_\alpha$ in the quotient group as
properties ``modulo $G_\alpha$'' of $X_\alpha$ in the group
$G_{\alpha+1}$ we obtain the following.

\begin{la}   \label{X-la}\mbox{}
  \begin{ls}
    \item $X_\alpha$ is linearly independent modulo $G_\alpha$.  That
    is, if $G_\alpha$ contains a linear combination, with integer
    coefficients, of members of $X_\alpha$, then all the coefficients
    are zero.
    \item $(\sq{X_\alpha}+G_\alpha)_*=G_{\alpha+1}$.
  \end{ls}
\end{la}

\begin{pf}
For the first assertion, notice that such a linear combination, when
projected to the quotient modulo $G_\alpha$, becomes a linear
combination of members of $Y_\alpha$ that equals zero.  So the linear
independence of $Y_\alpha$ in the quotient group gives the required
conclusion.

For the second assertion, consider an arbitrary element $a\in
G_{\alpha+1}$. Its image $\bar a$ in $G_{\alpha+1}/G_\alpha$ has a
multiple $n\bar a\in\sq{Y_\alpha}$ for some non-zero $n\in\bbb Z$.
Since the projection maps $X_\alpha$ onto $Y_\alpha$ and thus maps
\sq{X_\alpha} onto \sq{Y_\alpha}, we have an element
$z\in\sq{X_\alpha}$ projecting to $n\bar a$.  So $na$ and $z$ project
to the same element, which means $na=z+g$ for some $g\in G_\alpha$.
This equation establishes that $a$ is in the purification of
$\sq{X_\alpha}+G_\alpha$, as required.
\end{pf}

Temporarily fix an arbitrary stationary subset $T$ of $\Gamma(G)$.
Using $T$ we define a subgroup $L$ of $G$ by 
\[
L=\Big\langle\bigcup_{\beta\in T}X_\beta\Big\rangle_*.
\]
Our immediate objective, and indeed the main part of our argument for
Theorem~\ref{main} when \nf G is uncountable, is to show that $L$ is
not free.

For this purpose, we shall use Lemma~\ref{gamma-la} and show that
$\Gamma(L)$ is stationary.  (So we shall have not only that $L$ is not
free but that $\nf L=\lambda$.)  In fact, we shall prove more, namely
that $\Gamma(L)$ contains all of the stationary set $T$ except for
some non-stationary subset.

To begin the analysis of $\Gamma(L)$, we must choose a filtration of
$L$ to use in the definition of $\Gamma(L)$, and two natural choices
present themselves.  One is the restriction to $L$ of the filtration
we already have for $G$, i.e., \fil{G_\alpha\cap L:\alpha<\lambda}.
The other is obtained from the way $L$ is generated by the
$X_\alpha$'s; this filtration is \fil{L_\alpha:\alpha<\lambda}, where 
\[
L_\alpha=\Big\langle\bigcup_{\substack{\beta\in
    T\\\beta<\alpha}}X_\beta\Big\rangle_*. 
\]
We shall want to use each of these occasionally.  Fortunately, as
mentioned earlier, these two filtrations (like any two filtrations
of the same group) agree on a club.  Let $T_1$ be the intersection of
$T$ with such a club.  To show that $\Gamma(L)$ contains all but a
non-stationary part of $T$, it suffices, since $T-T_1$ is
non-stationary, to prove that $\Gamma(L)$ contains all but a
non-stationary part of $T_1$.

So our objective is now to show that the ``exceptional'' set 
\[
W=\{\alpha\in T_1:\alpha\notin\Gamma(L)\}
\]
is not stationary.  Suppose, toward a contradiction, that $W$ is
stationary.  

For each $\alpha\in W$, we have the following two facts:
\begin{ls}
  \item $L_{\alpha+1}/L_\alpha$ is free.
  \item $L_\alpha=L_{\alpha+1}\cap G_\alpha$.
\end{ls}
The first of these is immediate because $\alpha\notin\Gamma(L)$.  The
second follows from
\[
L_\alpha\subseteq L_{\alpha+1}\cap G_\alpha
\subseteq L\cap G_\alpha = L_\alpha,
\]
where the first inclusion uses the fact that $L_\alpha$ is the pure
subgroup of $G$ generated by a subset $\bigcup_{\beta\in T,\
  \beta<\alpha}X_\beta$ of $G_\alpha$ and $G_\alpha$ is a pure
subgroup of $G$.  The second inclusion is trivial, and the final
equality uses the fact that $\alpha\in T_1$ so the two filtrations
agree at $\alpha$.  

Combining the two facts just established, we have, for each $\alpha\in
W$, that the following group is free:
\[
\frac{L_{\alpha+1}}{L_{\alpha+1}\cap G_\alpha}\cong
\frac{L_{\alpha+1}+G_\alpha}{G_\alpha}.
\]
Contrast this with the fact that, since $\alpha\in
T\subseteq\Gamma(G)$, the following group is not free:
\[
\frac{G_{\alpha+1}}{G_\alpha}=
\frac{(\sq{X_\alpha}+G_\alpha)_*}{G_\alpha}.
\]
(Here we use that our filtration was chosen to satisfy the second
conclusion of Lemma~\ref{fil-norm}.)  Since subgroups of free groups
are free, we conclude that
\[
(\sq{X_\alpha}+G_\alpha)_*\not\subseteq L_{\alpha+1}+G_\alpha.
\]

Choose, for each $\alpha\in W$, some element
$g_\alpha\in(\sq{X_\alpha}+G_\alpha)_*$ that is not in
$L_{\alpha+1}+G_\alpha$.  By definition of purification, we have some
$n_\alpha\in\bbb Z-\{0\}$ such that
$n_\alpha g_\alpha\in\sq{X_\alpha}+G_\alpha$.  That is, $n_\alpha
g_\alpha = c_\alpha+h_\alpha$ where $c_\alpha$ is a linear
combination, with integer coefficients, of elements of $X_\alpha$ and
where $h_\alpha\in G_\alpha$.  

Because we assumed, toward a contradiction, that $W$ is stationary,
and because $h_\alpha\in G_\alpha$ for all $\alpha\in W$,
Lemma~\ref{multifodor} gives us a stationary set $W_1\subseteq W$ such
that $h_\alpha$ is the same element $h$ for all $\alpha\in W_1$.
Furthermore, since the values of $n_\alpha$ all lie in the countable
set \bbb Z, there is a stationary $W_2\subseteq W_1$ such that
$n_\alpha$ has the same value $n$ for all $\alpha\in W_2$.  

Let $\alpha<\beta$ be any two elements of $W_2$.  Then we have
\begin{align*}
  ng_\alpha &= c_\alpha+h\\ ng_\beta &= c_\beta+h.
\end{align*}
Subtract to get
\[
n(g_\beta-g_\alpha)=c_\beta-c_\alpha.
\]
Both $c_\alpha$ and $c_\beta$ are linear combinations of elements of
$X_\alpha\cup X_\beta\subseteq L_{\beta+1}$.  So $L_{\beta+1}$
contains the right side of the last equation.  As $L_{\beta+1}$ is
pure in $G$, it follows that $g_\beta-g_\alpha\in L_{\beta+1}$ and so
$g_\beta\in g_\alpha+L_{\beta+1}$.

Notice that $G_\beta$ includes $G_\alpha$ (because $\alpha<\beta$),
and $X_\alpha$ (because $X_\alpha\subseteq G_{\alpha+1}$ and
$\alpha+1\leq\beta$), and therefore $\sq{X_\alpha}+G_\alpha$, and
therefore $(\sq{X_\alpha}+G_\alpha)_*$ (because, according to the
first conclusion of Lemma~\ref{fil-norm}, $G_\beta$ is a pure subgroup
of $G$).  We chose $g_\alpha$ from this last group, so we have
$g_\alpha\in G_\beta$.  But then the result from the preceding
paragraph gives us that $g_\beta\in G_\beta+L_{\beta+1}$, which
contradicts our choice of $g_\beta$.

This contradiction completes the proof that $W$ cannot be stationary
and therefore $L$ is not free.  

Recall that $L$ was defined in terms of an arbitrary but fixed
stationary $T\subseteq\Gamma(G)$.  We shall now need to vary $T$, so,
to indicate the dependence of $L$ on $T$, we write $L(T)$ for what was
previously called simply $L$.

\begin{la}   \label{l-disjoint}
  If $T_1$ and $T_2$ are disjoint stationary subsets of $\Gamma(G)$,
  then the subgroups $L(T_1)$ and $L(T_2)$ are disjoint.
\end{la}

\begin{pf}
Notice first that, if two subgroups $H_1$ and $H_2$ are disjoint, then
so are their purifications.  Indeed, if the purifications had a common
non-zero element $x$, then $x$ would have non-zero multiples $n_1x\in
H_1$ and $n_2x\in H_2$.  But then $n_1n_2x$ would be a non-zero
element in $H_1\cap H_2$ contrary to hypothesis.

So to prove the lemma, it suffices to prove that the subgoups
generated by $\bigcup_{\beta\in T_1}X_\beta$ and $\bigcup_{\beta\in
T_2}X_\beta$ are disjoint.  Suppose, toward a contradiction, that they
are not disjoint, choose a non-zero element in their intersection, and
write it first as a linear combination $c_1$ of elements of
$\bigcup_{\beta\in T_1}X_\beta$ and second as a linear combination
$c_2$ of elements of $\bigcup_{\beta\in T_2}X_\beta$.  Let $\beta$ be
the largest of the finitely many ordinals such that elements of
$X_\beta$ occur in these linear combinations with non-zero
coefficients.  Since $T_1$ and $T_2$ are disjoint, $\beta$ is in only
one of them, say $T_1$.  Split $c_1$ as $c_1'+c_1''$ where $c_1'\neq0$
contains the terms from $X_\beta$ and $c_1''$ contains the terms from
$X_\alpha$'s with $\alpha<\beta$.  Then $c_1'=c_2-c_1''\in G_\beta$.
This contradicts the linear independence of $X_\beta$ modulo $G_\beta$
in Lemma~\ref{X-la}.
\end{pf}

Because of this lemma and the non-freeness of $L(T)$ for all
stationary $T\subseteq\Gamma(G)$, we can get as many pairwise
disjoint, non-free subgroups of $G$ as we can get pairwise disjoint,
stationary subsets of $\Gamma(G)$.  It remains only to quote Solovay's
famous theorem \cite[Theorem~9]{solovay} that, for any uncountable
regular cardinal $\lambda$, every stationary subset of $\lambda$ can
be partitioned into $\lambda$ pairwise disjoint, stationary subsets.

\section{Proof for Singular Non-Free Rank}  \label{sing}

In this section, we complete the proof of Theorem~\ref{main} by
treating the case of singular \nf G.  An important ingredient of the
proof is the following consequence of the second author's singular
compactness theorem.

\begin{la}   \label{sing-nfr}
Assume that $\kappa<\lambda$ are infinite cardinals, that $\lambda$ is
singular, that $G$ is a group of cardinality $\lambda$, and that every
subgroup $H$ of $G$ with $|H|<\lambda$ has $\nf H\leq\kappa$.  Then
$\nf G\leq\kappa$.
\end{la}

\begin{pf}
  This follows from the singular compactness theorem of
  \cite{sing-cpt}.  It is explicitly stated in \cite[Theorem~2]{hhs}.
  (To avoid possible confusion, note that in \cite{hhs}
  ``$\kappa$-generated'' means generated by strictly fewer than
  $\kappa$ elements.)  For an easier proof, see \cite{hodges}.
\end{pf}

Our use of this lemma will be via the following consequence.

\begin{la}   \label{scpt}
  Let $G$ be a group whose non-free rank is a singular cardinal
  $\lambda$.  For any infinite cardinal $\kappa<\lambda$, there exists
  a subgroup $H$ of $G$ whose non-free rank is regular and satisfies
  $\kappa<\nf H=|H|<\lambda$.
\end{la}

\begin{pf}
Given $G$, $\lambda$, and $\kappa$ as in the statement of the lemma,
let $H$ be a subgroup of $G$ with \nf H as small as possible subject
to the constraint that $\nf H>\kappa$.  Such an $H$ certainly exists,
since $\nf G>\kappa$.  By definition of \nf H, we can split $H$ as the
direct sum of a free group and a group $H_1$ of cardinality $\nf H$.
Then $|H_1|=\nf H=\nf{H_1}$.  Replacing $H$ with $H_1$ if necessary, we
assume from now on that $|H|=\nf H$.

It remains only to show that \nf H is a regular cardinal and that $\nf
H<\lambda$.  The latter follows from the former, because $\lambda$ is
singular, and so we need only prove the regularity of \nf H.  Suppose,
toward a contradiction, that \nf H is a singular cardinal.  Notice
that, for every subgroup $K$ of $H$ with cardinality smaller than
$|H|$, we have $\nf K\leq|K|<|H|=\nf H$, and so, by minimality of \nf
H, we must have $\nf K\leq\kappa$.  By Lemma~\ref{sing-nfr}, $\nf
H\leq\kappa$ also, but this contradicts the choice of $H$.
\end{pf}

Using the lemma, we prove the singular case of Theorem~\ref{main} as
follows.  Let $G$ be a group with $\nf G=\lambda$ singular, let
$\mu<\lambda$ be the cofinality of $\lambda$, and let
$\fil{\kappa_\xi:\xi<\mu}$ be a strictly increasing $\mu$-sequence of
cardinals with supremum $\lambda$.  Inductively define a $\mu$-sequence of
subgroups $G_\xi$ of $G$, with the following properties for all
$\xi<\mu$: 
\begin{ls}
  \item $|G_\xi|=\nf{G_\xi}$.
  \item \nf{G_\xi} is a regular cardinal.
  \item $\kappa_\xi<\nf{G_\xi}<\lambda$.
  \item $\sum_{\eta<\xi}|G_\eta|<\nf{G_\xi}$.
\end{ls}
Once $G_\eta$ is defined and has the desired properties for all
$\eta<\xi$, we obtain $G_\xi$ by applying Lemma~\ref{scpt} with
$\kappa$ equal to the larger of $\kappa_\xi$ and
$\sum_{\eta<\xi}|G_\eta|$.  Notice that this sum is strictly smaller
than $\lambda$ because $|G_\eta|<\lambda$ for all $\eta<\xi$ by
induction hypothesis and because $\xi<\mu=\text{cf}(\lambda)$.  Thus,
our $\kappa$ is smaller than $\lambda$, and the $H$ provided by
Lemma~\ref{scpt} serves as the required $G_\xi$.

The regular case of Theorem~\ref{main}, proved in the previous
section, applies to each $G_\xi$ and provides a family $\scr D_\xi$ of
\nf{G_\xi} pairwise disjoint non-free subgroups of $G_\xi$.  

Although the groups in any $\scr D_\xi$ are pairwise disjoint, there
may be a non-zero intersection between a group $H\in\scr D_\xi$ and a
group $K\in\scr D_\eta$ for some $\eta<\xi$.  But this does not happen
too often.  Specifically, for fixed $\eta$, $K$, and $\xi$, the number
of such $H$'s is at most $|K|\leq|G_\eta|$ because different $H$'s
from $\scr D_\xi$ would meet $K-\{0\}$ in disjoint sets.  Therefore,
if we keep $\eta$ and $\xi$ fixed but let $K$ vary through all
elements of $\scr D_\eta$, then the number of $H\in\scr D_\xi$ that
have non-zero intersection with some such $K$ is at most
$|G_\eta|\cdot|\scr D_\eta|$.  Since $\scr D_\eta$ is a family of
disjoint subgroups of $G_\eta$, its cardinality is at most $|G_\eta|$,
and so the product $|G_\eta|\cdot|\scr D_\eta|$ is also at most
$|G_\eta|$.  Now keeping $\xi$ fixed but letting $\eta$ vary through
all ordinals $<\xi$, we find that the number of $H\in\scr D_\xi$ that
have non-zero intersection with some $K$ in some earlier $\scr D_\eta$
is at most
\[
\sum_{\eta<\xi}|G_\eta|<\nf{G_\xi}.
\]
Therefore, if we discard these $H$'s from $\scr D_\xi$, what remains
is a family $\scr D'_\xi$, still of cardinality \nf{G_\xi}, still
consisting of pairwise disjoint, non-free subgroups of $G_\xi$, but
enjoying the additional property that all its members are disjoint
from all members of earlier $\scr D_\eta$'s and, a fortiori, from all
members of earlier $\scr D'_\eta$'s.

Therefore, $\scr D'=\bigcup_{\xi<\mu}\scr D'_\xi$ is a family of
pairwise disjoint non-free subgroups of $G$.  Its cardinality
satisfies 
\begin{align*}
|\scr D'|&=\sup\{|\scr D'_\xi|:\xi<\mu\}\\
&=\sup\{\nf{G_\xi}:\xi<\mu\}\\
&\geq\sup\{\kappa_\xi:\xi<\mu\}\\
&=\lambda,
\end{align*}
so the theorem is proved.

\section{Counterexample for Countable Non-Free Rank}  \label{counter} 

In this section, we construct an example of a group $G$ with $|G|=\nf
G=\aleph_0$ and with no two disjoint non-free subgroups.  This shows
that the uncountability assumption in Theorem~\ref{main} cannot be
removed. 

The group $G$ will be a subgroup of the direct sum of $\aleph_0$
copies of the additive group \bbb Q of rational numbers.  $G$ will be
defined as the set of solutions, in this direct sum, of infinitely
many $p$-adic conditions for all primes $p$.  The construction of $G$
will proceed in three phases.  The first will set up some conventions
and bookkeeping.  The second will define, for each prime $p$, a
certain set of vectors over $\bbb Z/p$ and a lifting of these vectors
to \zp.  The third will use these vectors to define the group $G$.
After the construction is complete, we shall verify first that $\nf
G=\aleph_0$ and second that $G$ does not have two disjoint, non-free
subgroups.  

\subsection{Conventions and bookkeeping}

By a \emph{vector} over a ring $R$, we shall usually mean an infinite
sequence of elements of $R$, the sequence being indexed by the set
\bbb N of positive integers.  (We do not include 0 in the index set,
as we shall have a separate, special use for 0 later.)  The rings
relevant to our work will include \bbb Z, \bbb Q, $\bbb Z/p$, \zp, and
\qp\ for prime numbers $p$.  

Occasionally, we shall need to refer to
vectors of finite length $l$, with components indexed by
$\{1,2,\dots,l\}$, but this finiteness (and the value of $l$) will
always be explicitly stated.  Furthermore, it would do no harm to
identify any finitely long vector with the infinite vector obtained by
appending a sequence of zeros.  We use the notation $\vec x\restr l$
for the vector of length $l$ consisting of just the first $l$
components of $\vec x$; under the identification in the preceding
sentence, $\vec x\restr l$ can also be considered as obtained from
$\vec x$ by replacing all components beyond the first $l$ by 0.  

We call a vector $\vec x$ \emph{finitely supported} if the set
$\{i\in\bbb N:x_i\neq0\}$, called the \emph{support} of $\vec x$, is
finite.  The inner product of two vectors is defined as
\[
\sq{\vec x,\vec y}=\sum_{i\in\bbb N}x_iy_i,
\]
provided at least one of the vectors is finitely supported, so that
the infinite sum makes sense.

Partition the set $P$ of prime numbers into infinitely many infinite
pieces, and label the pieces as $P_{\vec x}$, where $\vec x$ ranges
over the non-zero, finitely supported, infinite vectors over \bbb Z.
(The number of such vectors is $\aleph_0$, so the indexing makes
sense.)  Fix this indexed partition of $P$ for the rest of the proof.

Also fix an enumeration, as \fil{\vec\lambda_i:i\in\bbb N}, of all the
finitely supported vectors over \bbb Q.  

Recall that \bbb Q is canonically identified with a subfield of \qp,
and that under this identification \bbb Z becomes a subring of \zp.
This inclusion of \bbb Z in \zp\ induces an isomorphism between the
quotients modulo $p$, $\bbb Z/p\cong\zp/p\zp$.  We write $[x]_p$, or
just $[x]$ when $p$ is clear from the context, for the equivalence
class of $x$ modulo $p$; here $x$ is in \bbb Z or \zp, and $[x]$ is in
$\bbb Z/p$.  We refer to $x$ as a \emph{representative} of $[x]$.  We
use the same notation for vectors; $[\vec x]_p$ is obtained from $\vec
x$ by reducing all components modulo $p$.

\subsection{Useful sets of vectors}

For this subsection, let $p$ be a fixed prime.  Later, the work we do
here will be applied to all primes, but it is notationally and
conceptually easier to begin with just one $p$.  Let $\vec x$ be the
unique vector such that $p\in P_{\vec x}$.  Recall that this $\vec x$
is a non-zero, finitely supported vector over \bbb Z.  Thus, $[\vec
  x]$ is a finitely supported (but possibly zero) vector over $\bbb
Z/p$.  

We define a finite set $M$ of vectors over $\bbb Z/p$ as follows.
Choose an integer $l$ larger than $p$ and all elements of the support
of $\vec x$.  $M$ will be described as a set of finite vectors, in
$(\bbb Z/p)^l$, but we really mean the infinite vectors obtained by
appending a sequence of zeros.

If $[\vec x]=\vec0$, then $M$ consists of all the vectors in
$(\bbb Z/p)^l$.  If, on the other hand, $[\vec x]\neq\vec0$, then we
proceed as follows.  Call an index $i$ or the corresponding vector
$\vec\lambda_i$ (from the enumeration fixed above) \emph{relevant} if
$i<p-1$ and $p$ does not divide the denominator of any component of
$\vec\lambda_i$.  For each relevant $i$, the components of
$\vec\lambda_i$ are in \zp, so it makes sense to reduce them modulo
$p$, obtaining a vector $[\vec\lambda_i]$ over $\bbb Z/p$.  Choose a
non-zero $a\in\bbb Z/p$ that is distinct from
\sq{[-\vec\lambda_i],[\vec x]} for all relevant $i$.  The requirement,
in the definition of relevance, that $i<p-1$ means that at most $p-2$
indices are relevant, so a suitable $a$ exists.  Now let $M$ consist
of all those $\vec m\in(\bbb Z/p)^l$ such that $\sq{\vec m,[\vec
    x]}=a$.  

Before proceeding further, we summarize the properties of $M$ that we
shall need later.

\begin{la}   \label{m-la}
  The set $M\subseteq(\bbb Z/p)^l$ defined here has the following
  properties. 
  \begin{lsnum}
\item For each $k\leq l$, the truncations to $k$ of the vectors in $M$
  span $(\bbb Z/p)^k$.
\item The same holds for the truncations to $k$ of the translates
  $[\vec\lambda_i]+M$ for each relevant $\vec\lambda_i$.
\item As $\vec m$ varies through $M$, all the inner products \sq{\vec
  m,[\vec x]} have the same value $a$.
\item This value $a$ is non-zero and different from
  \sq{[-\vec\lambda_i],[\vec x]} for all relevant $i$ provided $[\vec
  x]\neq\vec0$. 
  \end{lsnum}
\end{la}

\begin{pf}
If $[\vec x]=\vec0$ then item~(3) is obvious with $a=0$ and item~(4)
doesn't apply.  If $[\vec x]\neq\vec0$ then items~(3) and (4) are
explicitly in the definition of $M$.  

As for items~(1) and (2), it suffices to prove these for $k=l$, since
truncation to smaller $k$'s commutes with linear combinations.  So
assume $k=l$.  If $[\vec x]=\vec0$ then both (1) and (2) are obvious
as $M$ is all of $(\bbb Z/p)^l$.  So assume $[\vec x]\neq\vec0$.  Then
$M$ is defined as a certain affine hyperplane in $(\bbb Z/p)^l$.
Since $a\neq0$, this hyperplane does not pass through the origin, and
therefore it spans $(\bbb Z/p)^l$.  This proves (1).  For (2), observe
that $[\vec\lambda_i]+M$ is another affine hyperplane, a translate of
$M$.  It does not pass through the origin either, because
$a\neq\sq{[-\vec\lambda_i],[\vec x]}$.  So it, too, spans $(\bbb
Z/p)^l$, and (2) is proved.
\end{pf}

The next step is to lift $M$ to an infinite set of vectors over \zp.
For this, we regard the vectors in $M$ as having infinite length
(rather than $l$), by appending zeros.  Arbitrarily choose, for each
$\vec m\in M$, a countable infinity of vectors $\vec\psi$ over \zp\
with $[\vec\psi]=\vec m$.  Let $\Psi$ be the set of all the resulting
vectors, for all $\vec m$ together.

Partition the set $\Psi$ into infinitely many pieces $\Psi_k$
($k\in\bbb N$) with $|\Psi_k|=k+1$.  For each $k$, we shall modify the
$k+1$ vectors in $\Psi_k$ so that their truncations to $k$ become
affinely independent over \qp.  Affine independence means that not
only do these $k+1$ vectors span $\qp^k$ but they continue to span if
any fixed vector is added to all of them.  Any set of $k+1$ vectors in
a $k$-dimensional vector space over \qp\ (or over any valued field
with a non-trivial valuation) can be made affinely independent by an
arbitrarily small perturbation.  For our modification of $\Psi_k$, we
use a perturbation that is small in the $p$-adic norm, so that each
vector $\vec\psi$ is modified by adding something divisible in \zp\ by
$p$.  Thus, the reduction modulo $p$, $[\vec\psi]$ is unaffected.  Let
the perturbed vectors constitute $\Phi_k$, and let
$\Phi=\bigcup_{k\in\bbb N}\Phi_k$.  Note that all the vectors in
$\Phi$ have, like those in $\Psi$, all their components in \zp.

\begin{la}  \label{phi-la}
  The set $\Phi$ of vectors over \zp\ constructed here has the
  following properties.
  \begin{lsnum}
\item For each $k\leq p$, the vectors $[\vec\phi]\restr k$ for
  $\vec\phi\in\Phi$ span $(\bbb Z/p)^k$.
\item The same holds for the translates $[\vec\lambda_i]+[\vec\phi]$,
  truncated at $k$, for each relevant $\vec\lambda_i$.
\item As $\vec\phi$ varies through $\Phi$, all the inner products
\sq{[\vec\phi],[\vec x]} have the same value $a$.
\item This value $a$ is non-zero and different from
  \sq{[-\vec\lambda_i],[\vec x]} for all relevant $i$ provided $[\vec
  x]\neq\vec0$.
\item For every $k$, there are $k+1$ vectors
  $\vec\phi_0,\vec\phi_1,\dots,\vec\phi_k\in\Phi$ such that, for any
  vector $\lambda$ over \bbb Q, the truncated vectors
  $(\vec\phi_i+\vec\lambda)\restr k$ span $\qp^k$.
  \end{lsnum}
\end{la}

\begin{pf}
  By construction, the reductions modulo $p$ of the vectors in $\Phi$
  are the same as those of the vectors in $\Psi$, namely the vectors
  in $M$ (each repeated $\aleph_0$ times).  Therefore, the first four
  conclusions of the present lemma follow immediately from the
  corresponding items in Lemma~\ref{m-la}.  (For the first item,
  remember that the $l$ in Lemma~\ref{m-la} was $\geq p$.)

The final conclusion of the present lemma is the result of our
modification of $\Psi$ to obtain $\Phi$.  The $k+1$ vectors in
$\Phi_k$ are affinely independent, and therefore their translates by
any vector $\vec\lambda$ over \qp\ are linearly independent.  This
applies in particular to vectors $\vec\lambda$ over $\bbb
Q\subseteq\qp$. 
\end{pf}

\subsection{The counterexample group}   

In the preceding subsection, we worked with a fixed prime $p$; now we
let $p$ vary.  The set $\Phi$ constructed above will now be called
$\Phi(p)$.    

Our group $G$ will be a subgroup of the additive group $\bbb
Q\times(\bbb Q)^{(\bbb N)}$, which consists of pairs $(x_0,\vec x)$
where $x_0$ is a rational number and $\vec x$ is a finitely supported
vector of rational numbers.  Recall that the components of
$\vec x$ are indexed as $x_i$ for positive integers $i$, so our use of
the notation $x_0$ causes no conflict.  We define
\[
G=\{(x_0,\vec x)\in\bbb Q\times(\bbb Q)^{(\bbb N)}:
(\forall p\in P)(\forall\vec\phi\in\Phi(p))\ 
x_0+\sq{\vec\phi,\vec x}\in\zp\}.
\]
Notice that the inner product in the definition makes sense because
$\vec x$ is finitely supported, the components of $\vec x$, being in
\bbb Q, can be regarded as elements of \qp, and the components of
$\vec\phi$, being in \zp, are also in \qp.  Thus, the inner product
is defined in \qp, and the requirement is that it be in \zp.

\begin{la}
  $\bbb Z\times(\bbb Z)^{(\bbb N)}\subseteq G$.
\end{la}

\begin{pf}
  Remember that all components of all $\vec\phi$ in $\Phi(p)$ are
  $p$-adic integers.  So if the components of $\vec x$ are integers,
  and thus also $p$-adic integers, then the inner product
  \sq{\vec\phi,\vec x} is a $p$-adic integer, and so is its
  sum with the integer $x_0$.  
\end{pf}

This lemma and the fact that $\bbb Q\times(\bbb Q)^{(\bbb N)}$ is
countable show that the rank and the cardinality of $G$ are
$\aleph_0$.  

We close this subsection by introducing a subgroup of $G$ that will
play a central role in our proof that $G$ has the desired properties.
Let $L$ be the subgroup of those elements of $G$ for which at most the
single component $x_0$ is non-zero.  So
\[
L=\{(x_0,\vec 0):(\forall p\in P)(\forall\vec\phi\in\Phi(p))\
x_0\in\zp\}. 
\] 
Recall that a rational number is in \zp\ for all primes $p$ if and
only if it is an integer.  Therefore, $L=\bbb Z\times\{0\}^{\bbb N}$.
Thus, $L$ is isomorphic to \bbb Z, and it is clearly a pure subgroup
of $G$.

\subsection{The non-free rank of $G$}

We now prove that $\nf G=\aleph_0$.  Since $|G|=\aleph_0$, we need
only show that $\nf G$ is infinite.  For this purpose, it is useful to
consider the quotient group $G/L$.

\begin{la}   \label{div-la}
  For each element $\xi\in G/L$, there are infinitely many primes $p$
  that divide $\xi$, i.e., such that $\xi\in p(G/L)$.
\end{la}

\begin{pf}
Given $\xi$, choose a representative $(x_0,\vec x)\in G$ for it modulo
$L$.  The components $x_i$ (both $x_0$ and the components of $\vec x$)
are rational numbers and only finitely many are non-zero, so let $d$
be a common denominator.  Then $d\xi$ is represented by a vector
$(dx_0,d\vec x)$ of integers.  If we prove the assertion of the lemma
for $d\xi$ then it will follow for $\xi$, using the same primes except
for the finitely many that divide $d$.  (In detail, if $d\xi=p\eta$
and $p$ doesn't divide $d$, then the Euclidean algorithm gives
integers $a,b$ with $ad+bp=1$.  Then $\xi=ad\xi+bp\xi=ap\eta+bp\xi$ is
divisible by $p$.)

So, by renaming, we may assume that all the components $x_i$ are
integers.  We may also assume that $\vec x\neq\vec0$, for
otherwise we would have $\xi=0$ in $G/L$ and the conclusion of the
lemma would be trivially correct.  So $P_{\vec x}$ is an infinite set
of primes.  We shall show that every $p\in P_{\vec x}$ divides $\xi$,
thereby completing the proof.

Since $\Phi(p)$ satisfies conclusion (3) of Lemma~\ref{phi-la}, we
know that the inner products \sq{\vec\phi,\vec x} have the same value
modulo $p$ for all $\phi\in\Phi(p)$.  Write this value as $[a]\in\bbb
Z/p$, where $a\in\bbb Z$, consider the element $(-a,\vec x)\in\bbb
Z\times(\bbb Z)^{(\bbb N)}\subseteq G$, and consider its quotient by
$p$,
\[
z=\left(-\frac ap,\frac 1p\vec x\right)\in\bbb Q\times(\bbb Q)^{(\bbb
  N)}. 
\]
Is this quotient $z$ in $G$?  Of the requirements for membership in $G$,
the $q$-adic ones for primes $q$ other than $p$ are certainly
satisfied, since division by $p$ preserves the property of being a
$q$-adic integer.  There remain the $p$-adic requirements.  But for
each $\phi\in\Phi(p)$, we have that $\sq{\vec\phi,\vec x}\equiv a\pmod p$
and therefore $-a+\sq{\vec\phi,\vec x}$ is divisible by $p$ in \zp.
This means that $-(a/p)+\sq{\vec\phi,(1/p)\vec x}\in\zp$, i.e., that
$z$ satisfies the $p$-adic
requirement arising from $\phi$.  Since this happens for every
$\phi\in\Phi(p)$, we conclude that $z\in G$.  But $pz$ differs from
$(x_0,\vec x)$ by an element of $L$, namely $(x_0+a,\vec 0)$.  So $pz$
represents $\xi$; that is, $\xi\in p(G/L)$ as required.
\end{pf}

\begin{cor}  \label{div-cor}
  $G/L$ has no non-zero, free, pure subgroup.
\end{cor}

\begin{pf}
  Any pure subgroup of $G/L$ would inherit from $G/L$ the divisibility
property in the lemma and therefore cannot be free unless it is the
zero group.
\end{pf}

\begin{la}   \label{nfrk-cble}
  The non-free rank of $G$ is $\aleph_0$. 
\end{la}

\begin{pf}
  As remarked above, we need only prove that \nf G cannot
  be finite.  Suppose, toward a contradiction, that $G$ can be split
  as $H\oplus F$, where $H$ has finite rank and $F$ is free.  Choose a
  free basis $B$ for $F$, and express the generator $(1,\vec 0)$ of
  $L$ as a linear combination of a vector from $H$ and a finite set
  $B_0$ of vectors from $B$.  Then $L\subseteq H\oplus\sq{B_0}$.  Let
  $B_1=B-B_0$ and notice that $B_1$ is infinite because $G$ has
  infinite rank, $H$ has finite rank, and $B_0$ is finite.  Since
  $G=H\oplus\sq{B_0}\oplus\sq{B_1}$, we have
\[
\frac GL\cong\left(\frac{H\oplus\sq{B_0}}L\right)\oplus\sq{B_1}.
\]
This makes \sq{B_1} a non-zero, pure, free subgroup of $G/L$, contrary
to Corollary~\ref{div-cor}.
\end{pf}

\subsection{Non-free subgroups of $G$}

To complete the verification of our counterexample, we must show that
$G$ does not have two disjoint, non-free subgroups.  The following
lemma is the main ingredient in that proof.

\begin{la}
If $H$ is a subgroup of $G$ of finite rank and $H\cap L=\{0\}$, then
$H$ is free.  
\end{la}

\begin{pf}
Let $H$ be as in the hypothesis.  Because it has finite rank and
because every element of $G$ has finite support, $H$ is a subgroup of
$\bbb Q\times\bbb Q^k\times\{\vec0\}$ for some $k$.  Whenever
convenient, we shall ignore the $\{\vec0\}$ factor and pretend
that elements of $H$ are of the form $(x_0,\vec x)$ with $\vec
x\in\bbb Q^k$.  

Recall (for example from the proof of Lemma~\ref{l-disjoint}) that
purifications of disjoint subgroups are disjoint.  Since $H$ is
disjoint from $L$, the purification $H_*$ of $H$ in $\bbb Q\times\bbb
Q^k$ is disjoint from the purification $L_*=\bbb Q\times\{0\}^k$ of
$L$.  These purifications are \bbb Q-linear subspaces of $\bbb
Q\times\bbb Q^k$, and so we know, by linear algebra, that there is a
vector $\vec\lambda\in\bbb Q^k$ such that all elements $(x_0,\vec x)$
of $H$ satisfy $x_0=\sq{\vec\lambda,\vec x}$.  Being in $G$, these
elements of $H$ also satisfy
\begin{eqnum}   \label{lam-la}
\sq{\vec\lambda+\vec\phi,\vec x}=x_0+\sq{\vec\phi,\vec x}\in\zp
\end{eqnum}
for all primes $p$ and all $\phi\in\Phi(p)$.  

By appending zeros, we regard $\vec\lambda$ as an infinite vector.
Being a vector over \bbb Q with finite support, it occurs in our
enumeration as $\lambda_i$ for a certain $i$, which we fix for the
rest of this proof.  

Call a prime $p$ \emph{good} if it has the following properties:
\begin{ls}
\item $p\geq k$.
\item $p-1>i$.
\item $p$ does not divide the denominator of any component of
  $\lambda=\lambda_i$.  
\end{ls}
We note for future reference that these conditions are satisfied by
all but finitely many primes.  We also note that the second and third
of these conditions make $\lambda$ relevant for $p$ in the sense used
in Lemmas~\ref{m-la} and \ref{phi-la}.

Temporarily fix a good prime $p$.  Assume, toward a contradiction,
that $p$ divides the denominator of some component of $\vec x$ for
some element $(x_0,\vec x)$ of $H$.  Fix such an $(x_0,\vec x)$, let
$p^m$ be the highest power of $p$ that divides the denominator of a
component of $\vec x$, and let $\vec y=p^m\vec x$.  Then all
components of $\vec y$ have denominators prime to $p$, and at least
one of these components also has its numerator prime to $p$.  Since
$(x_0,\vec x)\in H$ and since $m>0$, formula~\eqref{lam-la} above
gives us, for every $\vec\phi\in\Phi(p)$,
\[
\sq{\vec\lambda+\vec\phi,\vec y}=
p^m\sq{\vec\lambda+\vec\phi,\vec x}\in p\zp,
\]
and so $\sq{[\vec\lambda+\vec\phi],[\vec y]}=0$ in $\bbb Z/p$.  Since
this holds for all $\vec\phi\in\Phi(p)$ and since, by conclusion~(2)
of Lemma~\ref{phi-la} the corresponding vectors
$[\vec\lambda+\vec\phi]\restr k$ span $(\bbb Z/p)^k$, it follows that
$[\vec y]\restr k$ is the zero vector in $(\bbb Z/p)^k$.  The
components of $\vec y$ beyond the first $k$ all vanish by our choice
of $k$.  Therefore, all components of $\vec y$ are divisible by $p$ as
$p$-adic integers, which means that, as rational numbers, they have 
numerators divisible by $p$.  That contradicts our choice of $m$
and $\vec y$.

This contradiction shows that good primes $p$ cannot divide
denominators of components of $\vec x$ when $(x_0,\vec x)\in H$.  Nor
can they divide $x_0=\sq{\vec\lambda,\vec x}$ because, being good,
they do not divide denominators of components of $\vec\lambda$.  Thus,
good primes do not divide the denominators of any components of
elements of $H$.  

We now turn our attention to those finitely many primes $p$ that are
not good.  Although it is possible for such a prime to divide the
denominator of a component of an element of $H$, we intend to show
that this divisibility cannot be with great multiplicity.  That is,
there is a bound $m\in\bbb N$ such that no power of $p$ higher than
$p^m$ divides the denominator of any element of $H$.

Fix one of these bad primes $p$, and choose $k+1$ vectors
$\vec\phi_j\in\Phi(p)$ as in assertion~(5) of Lemma~\ref{phi-la}.
Applying that assertion with our present $\vec\lambda$, we obtain that
the $k+1$ vectors $(\vec\phi_j+\vec\lambda)\restr k$ span ${\qp}^k$.
As ${\qp}^k$ is $k$-dimensional, we can select $k$ of these vectors
$\zeta_j=\vec\phi_j+\vec\lambda$ whose truncations form a basis for
${\qp}^k$.  This means that the matrix $Z$ whose rows are these
truncations is a non-singular $k\times k$ matrix over \qp.  (We use
here that the components of $\vec\phi_j$ are in $\zp\subseteq\qp$ and
the components of $\vec\lambda$ are in $\bbb Q\subseteq\qp$.)  Recall
that any $p$-adic number can be written as a $p$-adic integer divided
by a power of $p$.  So we can choose an integer $m$ so large that all
entries of the matrix $p^mZ^{-1}$ are in \zp.

For any $(x_0,\vec x)\in H\subseteq G$, we have, by
formula~\eqref{lam-la}, that each $\sq{\vec\zeta_j,\vec x}\in\zp$,
which means that, if we regard $\vec x$ as a column vector, then
$Z\vec x$ is a vector $\vec z$ with components in \zp.  Then $\vec
x=Z^{-1}\vec z$ has all its components in $p^{-m}\zp$.  That is, the
denominator of any component of $\vec x$ cannot be divisible by a
higher power of $p$ than $p^m$.  

To get the same result for $x_0$, we may need to increase $m$ as
follows.  Since $\vec\lambda$ has finite support, let $p^r$ be the
highest power of $p$ that divides the denominator of any component of
$\vec\lambda$.  Then, since $x_0=\sq{\vec\lambda,\vec x}$, no power of
$p$ higher than $p^{m+r}$ can divide the denominator of $x_0$.

Summarizing, we have an upper bound for the powers of $p$ that can
divide the denominator of any component of any member of $H$.

Let $D$ be the product of these powers $p^{m+r}$ for all of the
finitely many bad primes.  What we have shown is that every component
of every element of $H$ has, as its denominator, a divisor of $D$.
Indeed, such a denominator cannot have good prime factors, and the
remaining primes, the bad ones, cannot divide such a denominator to a
higher power than they divide $D$.  This means that $H$ is a subgroup
of
\[
\frac 1D\left(\bbb Z\times\bbb Z^k\times\{\vec 0\}\right).
\]
This group, isomorphic to $\bbb Z^{k+1}$, is free, and therefore so is
$H$. 
\end{pf}

\begin{cor}
  Every subgroup of $G$ that is disjoint from $L$ is free.
\end{cor}

\begin{pf}
  If $H$ is such a subgroup, then all its finite-rank subgroups are
  free, by the lemma.  Also, $H$ is countable, because $G$ is.  By
  Pontryagin's criterion (\cite[Lemma~16]{pont} or
  \cite[Theorem~IV.2.3]{em}), it follows that $H$ is free.
\end{pf}

Finally, we deduce that $G$ does not have two disjoint non-free
subgroups. Let two non-free subgroups of $G$ be given.  By the
corollary just proved, each of them contains a non-zero element of
$L$.  As $L\cong\bbb Z$, these elements have a common non-zero
multiple, which has to be in both of the given subgroups.

\end{document}